\documentclass[A4]{amsart}

\usepackage[all]{xy}
\usepackage{amssymb}  
\usepackage{latexsym} 
\usepackage{amsmath}
\usepackage{amsthm}

\DeclareFontEncoding{OT2}{}{} % to enable usage of cyrillic fonts
\newcommand{\textcyr}[1]{%
 {\fontencoding{OT2}\fontfamily{wncyr}\fontseries{m}\fontshape{n}\selectfont #1}}

\newcommand{\Sha}{{\mbox{\textcyr{Sh}}}}

\newcommand{\Z}{{\mathbb Z}}
\newcommand{\Q}{{\mathbb Q}}
\newcommand{\R}{{\mathbb R}}
\newcommand{\F}{{\mathbb F}}

\newcommand{\G}{{\mathbb G}}

\newcommand{\To}{\longrightarrow}
\newcommand{\rk}{\operatorname{rank}}

\newcommand{\ord}{\operatorname{ord}}

\newcommand{\GL}{\operatorname{GL}}

\newcommand{\HH}{\operatorname{H}}

\newcommand{\Hom}{\operatorname{Hom}}

\newcommand{\Pic}{\operatorname{Pic}}

\newcommand{\res}{\operatorname{res}}

\newcommand{\Br}{\operatorname{Br}}
\newcommand{\inv}{\operatorname{inv}}

\newenvironment{Proof}{\par\noindent{\sc Proof:}}%
                      {\hspace*{\fill}\nobreak$\Box$\par\hspace{2mm}}
\newenvironment{Remark}{\par\noindent{\sc Remark:}}{\par\hspace{2mm}}
                      {\hspace*{\fill}\nobreak$\Box$\par}

\newtheorem{Theorem}{Theorem}[section]
\newtheorem{Lemma}[Theorem]{Lemma}
\newtheorem{Proposition}[Theorem]{Proposition}
\newtheorem{Corollary}[Theorem]{Corollary}

\newtheorem{Question}[Theorem]{Question}

\begin{document}

\title{Potential Sha for abelian varieties}
\author{Brendan Creutz}
\address{School of Mathematics and Statistics, University of Sydney, NSW 2006, Australia}

\begin{abstract}
We show that the $p$-torsion in the Tate-Shafarevich group of any principally polarized abelian variety over a number field is unbounded as one ranges over extensions of degree $\mathcal{O}(p)$, the implied constant depending only on the dimension of the abelian variety.
\end{abstract}

\maketitle

\section{Introduction}
Let $A/k$ be an abelian variety over a number field $k$. We denote the completion of $k$ at a prime $v$ by $k_v$. The Galois cohomology group $\HH^1(k,A)$ is known as the Weil-Ch\^atelet group. The elements in this group can be interpreted as isomorphism classes of $k$-torsors under $A$. Such a torsor is trivial if and only if it contains a $k$-rational point. Since any torsor under the connected group $A$ is locally solvable at all but finitely many primes, the image of the restriction map 
\[ \res:\HH^1(k,A) \to \prod_v \HH^1(k_v,A) \] is contained in the direct sum. The Tate-Shafarevich group of $A/k$, denoted $\Sha(A/k)$, is defined as the kernel of this restriction map. As such, a torsor represents an element of the Tate-Shafarevich group if and only if it is everywhere locally solvable. Thus, nontrivial elements in this group are counterexamples to the Hasse principle.

Conjecturally the Tate-Shafarevich group is finite. Nevertheless there are a number of results in the literature which show that this group can be arbitrarily large. The first results in this direction were due to Cassels \cite{Ca1} who showed that one can arrange for the $3$-torsion subgroup of $\Sha(E/\Q)$ to be as large as one likes by choosing a suitable CM elliptic curve $E$ defined over $\Q$. Similar results have been obtained for all $p \le 7$ and $p = 13$ (see \cite{Boelling, Donnelly, Fisher, Kramer, Matsuno}). These results exploit the existence of elliptic curves admitting $p$-isogenies and, for $p = 5,7,13$, make use of the fact that the modular curve $X_0(p)$ parametrizing such curves has genus $0$ and infinitely many $\Q$-points. Kloosterman and Schaefer \cite{Kloosterman, KlSc} have employed similar techniques when the modular curve $X_0(p)$ is of higher genus to show that the order of the $p$-torsion in $\Sha(E/k)$ is unbounded as $E/k$ ranges over elliptic curves defined over number fields of degree bounded by some polynomial in $p$.

These latter results require that both the elliptic curve and the base field be allowed to vary. Matsuno showed \cite{Matsuno2} that if $k$ is a cyclic extension of $\Q$ of degree $p$, then the $p$-rank of the Tate-Shafarevich group over $k$ of elliptic curves defined over $\Q$ can be arbitrarily large. In a different direction, results of Clark and Sharif explore the $p$-torsion in $\Sha(A/\ell)$ for a fixed abelian variety $A$ over a number field $k$ as the extension $\ell/k$ is allowed to vary. Their results actually concern the subgroup  \[ \Sha_k(A/\ell) := \res_{\ell/k}(\HH^1(k,A))\cap\Sha(A/\ell)\,,\] which we call the {\em potential $\Sha$ of $A/k$ in $\ell$}. 

In \cite{WC1} Clark showed that the potential $p$-torsion in $\Sha$ of an elliptic curve $E/k$ with full level $p$ structure is unbounded as $\ell$ ranges over extensions of degree $p$. In \cite{ClSh} Clark and Sharif removed the assumption of full level $p$ structure, thus showing that the potential $p$-torsion of $\Sha$ of any elliptic curve is unbounded as $\ell$ ranges over extensions of degree $p$. Clark has generalized the original result (i.e. under the assumption of full level $p$ structure) to all strongly principally polarized abelian varieties for which the absolute Galois group $G_k$ of $k$ acts trivially on the N\'eron-Severi group \cite[Theorem 9]{ClarkAV} (recall that a polarization is {\em strong} if it is given by a $k$-rational divisor). The main result of this paper is to remove the assumption of full level $p$ structure in the higher-dimensional case as well.

\begin{Theorem}
\label{potentialSha}
Let $A/k$ be a strongly principally polarized abelian variety over a number field $k$ such  that the $G_k$-action on the N\'eron-Severi group is trivial. For any prime $p$ and any integer $N$, there exists a degree $p$ extension $\ell$ of $k$ for which \[ \#\Sha(A/\ell)[p] \ge \#\Sha_k(A/\ell)[p] > N\,.\]
\end{Theorem}

As Clark has noted (see \cite[Corollary 10]{ClarkAV}), any principally polarized abelian variety of dimension $g$ can be made to satisfy the conditions of theorem \ref{potentialSha} by passing to an extension of degree at most $2^g\cdot\#\GL_{4g^2}(\F_3)$. As this does not depend on $p$, we see that the $p$-torsion in the Tate-Shafarevich group of any principally polarized abelian variety $A/k$ becomes arbitrarily large as $\ell$ ranges over extensions of degree $\mathcal{O}(p)$, the implied constant depending only on the dimension of the abelian variety.

Our approach to theorem \ref{potentialSha} is based on the method utilized in \cite{WC1, ClarkAV, ClSh} which relates potential $\Sha$ to the period-index problem for torsors under abelian varieties. Recall that the {\em period} of a torsor under $A/k$ is its order in the Weil-Ch\^atelet group while its {\em index} is the greatest common divisor of the degrees of extensions over which the torsor contains a rational point. Using weak approximation together with a result of Lang and Tate regarding period and index over local fields (theorem \ref{LTtheorem}) one can often arrange that a torsor of degree $n$ become everywhere locally solvable after a sufficiently ramified extension $\ell/k$ of degree $n$. If the index of such a torsor exceeds $n$, then it will represent a nontrivial element of $\Sha(A/\ell)[n]$. Of course the difficulty lies in constructing sufficiently many such torsors and ensuring that they remain pairwise nonisomorphic over $\ell$.  We give a new construction of such period-index discrepancies which simultaneously generalizes \cite[Theorem 8]{ClarkAV} and \cite[Theorem 2]{ClSh}.

\begin{Theorem}
\label{Thm1}
Let $n \ge 2$ and $A/k$ a strongly principally polarized abelian variety over a number field $k$ such that the $G_k$-action on the N\'eron-Severi group is trivial. Then there exist infinitely many elements of $\HH^1(k,A)$ of period $n$ which cannot be split by any extension of degree $n$. If $n$ is a prime power, then these torsors have index strictly larger than $n$.
\end{Theorem}

For the statement of the next two results, let $A/k$ be a principally polarized abelian variety over a number field $k$, and let $p$ be a prime number. Let $S$ denote the set of primes of $k$ that are either of bad reduction for $A$ or lie above $p$, and let $T = \{ v \notin S \,:\, A(k_v)[p] \ne 0 \}$. Matsuno showed \cite[Corollary 4.5]{Matsuno2} that for an elliptic curve $E/\Q$ and a cyclic extension $\ell/\Q$ of degree $p$, \[ \#\Sha(E/\ell)[p] \cdot \#\left(\frac{E(\ell)}{pE(\ell)}\right) \ge p^{t - 4}\,,\] where $t$ denotes the number of primes in $T$ which are totally ramified in $\ell$. Our approach yields similar results on the potential $p$-torsion in terms of the ramification of the extension.

\begin{Theorem}
\label{bounds1}
Assume $A$ is strongly principally polarized and that the $G_k$-action on the N\'eron-Severi group is trivial. Let $\{\ell_i\}$ be any sequence of degree $p$ extensions of $k$ such that for every prime $v$ in $T$, there exists an integer $n$ such that for every $i \ge n$, $v$ is totally ramified in $\ell_i$. Then $\lim_{i \to \infty} \#\Sha_k(A/\ell_i)[p]= \infty$
\end{Theorem}

\begin{Theorem}
\label{bounds2}
Let $\ell/k$ be any extension of degree $p$. Let $t$ be the number of primes in $T$ which are totally ramified in $\ell$. Then \[\#\Sha_k(A/\ell)[p^\infty] \le \#\Sha(A/k)[p^\infty] \cdot p^{g[k:\Q]+2g(t+\#S)}\,.\]
\end{Theorem}

\begin{Remark}
The reader will note that theorem \ref{bounds2} makes no assumptions on either the principal polarization or the $G_k$-module structure of the N\'eron-Severi group. The proofs of theorems \ref{potentialSha}--\ref{bounds1} rely on a theorem of Clark (\ref{ObThm} below) whose proof does require such assumptions.
\end{Remark}

The existence of an upper bound as in theorem \ref{bounds2} was alluded to by Clark and Sharif in \cite[Section 3.8]{ClSh}. They also ask if the size of $\Sha(E/\ell)[p]$ necessarily approaches infinity with the number of primes ramifying in $\ell$. Theorems \ref{bounds1} and \ref{bounds2} show that for the potential $\Sha$ it is only the ramification in $T$ that plays any role. In a sense this is quite natural. One can show that every torsor $V/k$ of $p$-powered period is locally solvable at all primes outside $S\cup T$.

\section{The obstruction map for abelian varieties}
Here we collect various results related to the {\em obstruction map} utilized by Clark to study the period-index problem for abelian varieties in \cite{ClarkAV}.  We are indebted to him for having laid a comfortable foundation for forays in this direction. His results generalize those of O'Neil who developed the obstruction map for elliptic curves \cite{ONeil}. The only possibly new material here is contained in propositions \ref{ObisTate} and \ref{ImageOb}. The first relates the obstruction map to the Tate pairing. The second uses this to compute the image of the obstruction map over $p$-adic fields. These results were to be expected, if not already well known.

Let $(A,L)$ be a polarized abelian variety over a field $K$. Recall that the polarization is given by a line bundle $L \in \Pic(\bar{A})$ that is ample, base point free and algebraically equivalent to each of its Galois conjugates. The last condition means that $L$ gives rise to a $G_K$-invariant class in the N\'eron-Severi group. One says that the polarization is {\em strong} if $L \in \Pic(A)$ (i.e. if $L$ can be represented by a $K$-rational divisor on $A$). One says that $L$ is {\em symmetric} if $L \simeq [-1]^*L$. Corresponding to $L$ is an isogeny $\varphi_L: A \to A^\vee$, given by $\varphi_L(x) = \tau_x^*L\otimes L^{-1}$, where $\tau_x$ denotes translation by $x$ on $A$. We denote the kernel of $\varphi_L$ by $A[\varphi_L]$ and let $n$ be the exponent of $A[\varphi_L]$. We assume that $n$ is not divisible by the characteristic of $K$.

\subsection{The result of Zarkhin}  Let $\mathcal{G}_L$ denote the theta group associated to $L$ (e.g. \cite{Mumford}, \cite[Section 5.1]{ClarkAV}, \cite[Chapter 8]{vanderGeer}). This is a central extension of $A[\varphi_L]$,
\begin{align}
\label{ThetaGroup}
1 \To \G_m \To \mathcal{G}_L \To A[\varphi_L] \To 0\,,
\end{align} which gives rise to a nondegenerate symplectic form \begin{align}
\label{SymplecticForm}
e^L:A[\varphi_L]\times A[\varphi_L] \To \G_m\,.
\end{align} For $Q,R \in A[\varphi_L]$ the pairing is given by taking the commutator of any choice of lifts of $Q$ and $R$ to $\mathcal{G}_L$. Since the possible lifts differ only by central elements the pairing is well-defined. Since $A[\varphi_L]$ is commutative, the pairing takes values in the center of $\mathcal{G}_L$. The fact that $\G_m$ is the center of $\mathcal{G}_L$ implies that the pairing is nondegenerate.

The pairing may be composed with the cup product to obtain a pairing \[ \cup_{e^L}: \HH^1(K,A[\varphi_L])\times\HH^1(K,A[\varphi_L]) \stackrel{\cup}{\To} \HH^2(K,A[\varphi_L]\otimes A[\varphi_L]) \stackrel{e^L}{\To} \HH^2(K,\G_m)\,.\] On the other hand, the central extension (\ref{ThetaGroup}) leads to an exact sequence (of pointed sets) coming from nonabelian Galois cohomology. The following relation between the connecting map \begin{align}
\label{DefDelta}
\Delta_L:\HH^1(K,A[\varphi_L]) \to \HH^2(K,\G_m) = \Br(K)
\end{align} and the cup product pairing was established by Zarkhin in \cite{Zarkhin}.

\begin{Theorem}[Zarkhin]
For any $\xi,\eta \in \HH^1(K,A[\varphi_L])$ one has \[ \xi \cup_{e^L} \eta = \Delta_L(\xi+\eta)-\Delta_L(\xi)-\Delta_L(\eta)\,.\] If $L$ is symmetric, then in addition $\Delta_L(a\xi) = a^2\Delta_L(\xi)$ for any integer $a$. In particular, $\Delta_L$ is quadratic and is a quadratic form if $L$ is symmetric.
\end{Theorem}

We record here the following elementary lemma for later use.

\begin{Lemma}
\label{QuadraticProperties}
Let $m$ be an odd number, $V$ a $\Z/m\Z$-module and $Q : V \to \Z/m\Z$ a quadratic form whose associated bilinear form $B$ is symmetric and non-degenerate. Suppose there exists an isotropic element $u \in V$ of order $m$ (i.e. such that $u \ne 0$ and $2Q(u) = B(u,u) = 0$). Then $Q(V) = \Z/m\Z$.
\end{Lemma}

\begin{Proof}
Let $a \in \Z/m\Z$. We will show that $a \in Q(V)$. Let $u$ be isotropic and of order $m$. The map $B(u,\cdot):V \to \Z/m\Z$ is surjective (since $u$ has order $m$ and $B$ is nondegenerate). Thus we can find $w$ such that $B(u,w) = 1/2 \in \Z/m\Z$. Set $b = a - B(w,w)$. Consider $v = bu+w$. We have
\begin{align*}
2Q(v) &= B(bu+w,bu+w)\\
&= 2\left(2bB(u,w)+B(w,w)+b^2B(u,u)\right)\\
&= 2\left(a - B(w,w) + B(w,w)\right) = 2a\,.
\end{align*} Since $2$ is invertible we have $Q(v) = a$.
\end{Proof}

\subsection{The Weil pairing} 
Let $\lambda:A \to B$ be an isogeny of abelian varieties. One has the dual isogeny $\lambda^\vee:B^\vee \to A^\vee$ and standard duality theorems (e.g. \cite[Theorem 7.5]{vanderGeer}) give an isomorphism of group schemes $\beta: B^\vee[\lambda^\vee] \simeq \Hom(A[\lambda],\G_m)$. The rule $(x,y) \mapsto \beta(y)(x)$ defines a pairing 
\begin{align}
\label{Weil1}
e_{\lambda}:A[\lambda]\times B^\vee[\lambda^\vee] \to \G_m\,. 
\end{align}
In the particular case $\lambda = n : A \to A$, this gives a pairing
\begin{align}
\label{WeilPairing}
e_n: A[n]\times A^\vee[n] \to \mu_n\,,
\end{align} known as the {\em Weil pairing}. In the situation considered above, one has an isogeny $\varphi_L:A \to A^\vee$. So $A[\varphi_L] = A[\varphi_L^\vee]$, and this gives rise to a pairing
\begin{align}
\label{Weil2}
e_{\varphi_L}:A[\varphi_L]\times A[\varphi_L] \to \G_m \,.
\end{align} The following proposition gives the connection between this and the commutator pairing in (\ref{SymplecticForm}).

\begin{Proposition}
\label{comparepairings}
$e^L = e_{\varphi_L}$.
\end{Proposition}

\begin{Proof}
This appears to be well known. It is alluded to in \cite{Zarkhin} and proven in \cite[Theorem 11.20]{vanderGeer}.
\end{Proof}

\subsection{The Tate pairing}
Together with the cup product, the Weil pairing $e_n$ induces a symmetric bilinear pairing
\begin{align}
\label{Tate1}
T_n:\HH^1(K,A[n])\times\HH^1(K,A^\vee[n]) \to \HH^2(K,\mu_n)=\Br(K)[n]\,.
\end{align} Recall that for any $n$ indivisible by the characteristic of $K$, we have a Kummer sequence
\begin{align}
\label{Kummer}
0 \To A(K)/nA(K) \stackrel{\delta}{\To} \HH^1(K,A[n]) \To \HH^1(K,A)[n] \To 0\,.
\end{align}
If $V \in \HH^1(K,A)[n]$ is a torsor under $A$ of period dividing $n$, then any lift to $\HH^1(K,A[n])$ will be called a {\em Kummer lift of $V$}. Tate has shown (e.g. \cite[Section I.3]{ADT}) that $T_n$ descends to give a symmetric bilinear pairing
\begin{align}
\label{Tate2}
\tilde{T}_n:\frac{A(K)}{nA(K)}\times\HH^1(K,A^\vee)[n] \to \Br(K)[n]\,,
\end{align} which is compatible with $T_n$ and the Kummer sequences of $A$ and $A^\vee$. Namely $\tilde{T}_n(R,V) = T_n(\delta(R), \tilde{V})$ where $\delta(R)$ is the image of $R$ in $\HH^1(K,A[n])$ under the connecting homomorphism and $\tilde{V}$ is any Kummer lift of $V$. 

\subsection{The case of a principal polarization}
Suppose now that $(A,P)$ is a principally polarized abelian variety over $K$, and consider the isogeny $\varphi_n$ associated to the line bundle $L = nP$. Upon identification of $A^\vee[n]$ and $A[n]$ using the principal polarization, the isogeny $\varphi_n$ is multiplication by $n$. Using the same identification, we can interpret the pairings $T_n$ and $\tilde{T}_n$ as being defined on $\HH^1(K,A[n])\times\HH^1(K,A[n])$ and $A(K)/nA(K)\times \HH^1(K,A)[n]$, respectively. For simplicity we denote the corresponding connecting map $\Delta_{L}$ in (\ref{DefDelta}) by $\Delta_n$.

\begin{Proposition}
\label{ObisTate}
Let $A/K$ be a principally polarized abelian variety over a field $K$ of characteristic not dividing $n$. Then the bilinear form associated to $\Delta_n$ is $T_n$.
\end{Proposition}

\begin{Proof}
By Zarkhin's result we know that the bilinear form associated to $\Delta_n$ is the cup product associated to the commutator pairing $e^{nP}$ for the theta group of $A[n]$. By proposition \ref{comparepairings} this coincides with the pairing $e_{\varphi_{nP}}$ associated to the isogeny $\varphi_{nP}$. As $\varphi_{nP}$ is multiplication by $n$ on $A \simeq A^\vee$, we see that $e_{\varphi_{nP}}$ is the Weil pairing $e_n:A[n]\times A[n] \to \mu_n$. The cup product induced by this is equal to $T_n$.
\end{Proof}

For a principally polarized abelian variety, the map
\[ \Delta_n : \HH^1(K,A[n]) \To \Br(K)\,,\] will be referred to as the {\em obstruction map}. Clark gives two other equivalent definitions of the obstruction map which allow him to use $\Delta_n$ to study the period-index problem for torsors under $A$. Rather than giving details, we encourage the reader to consult \cite[Section 5]{ClarkAV}. The following theorem gives the crucial connection between this map and the index of a torsor.

\begin{Theorem}[Clark] 
\label{ObThm}
Let $A/K$ be a strongly principally polarized abelian variety over a field $K$ of characteristic not dividing $n$. Assume $A$ has $G_K$-invariant N\'eron-Severi group and let $V \in \HH^1(K,A)[n]$. If $\Delta_n(\xi)\ne 0$ for every Kummer lift $\xi$ of $V$, then $V$ cannot be split over any degree $n$ field extension. If $n$ is a prime power, we have moreover that the index of $V$ exceeds $n$.
\end{Theorem}
\begin{Proof} 
See \cite[Theorem 31]{ClarkAV}. The proof there uses the assumption that the principal polarization is strong and that the N\'eron-Severi group is trivial as a Galois module.
\end{Proof}

\subsection{Local results}
\label{LocalResults}
We now specialize to the case that $A/k$ is a principally polarized abelian variety defined over a number field $k$. The definition of the obstruction map is functorial in the base field. In particular it commutes with restriction maps. For a completion $k_v$ of $k$ at a prime $v$, we we will use $\Delta_{n,v}$ to denote the corresponding `local obstruction map' $\Delta_{n,v}:\HH^1(k_v,A[n]) \to \Br(k_v)$. Similarly the local versions of the Tate pairings $T_n$ and $\tilde{T}_n$ will be denoted by $T_{n,v}$ and $\tilde{T}_{n,v}$, respectively.

We recall the following important theorem of Tate (see \cite[I.3.4]{ADT}).

\begin{Theorem}[Tate]
\label{TateLocal}
For any nonarchimedean prime $v$, the pairings $T_{n,v}$ and $\tilde{T}_{n,v}$ are nondegenerate.
\end{Theorem}

The Brauer group of $k_v$ embeds in $\Q/\Z$. Thus these pairings establish Pontryagin duality between the finite groups appearing in their domains. In the nonarchimedean case, one can determine the size of these groups by using the fact that $A(k_v)$ contains a finite index subgroup isomorphic to $\dim(A)$ copies of the maximal ideal in the ring of integers of $k_v$.

\begin{Lemma}
\label{sizetorsors}
Let $v$ be a nonarchimedean prime of $k$ above the rational prime $q$, and let $n$ be a positive integer. Then
\[ \#\HH^1(k_v,A)[n] = \#A(k_v)/nA(k_v) = q^{\dim(A)\cdot\ord_q(n)\cdot[k_v:\Q_q]} \cdot\#A(k_v)[n]\,.\]
\end{Lemma}

\begin{Proof}
The first equality follows from the nondegeneracy of $\tilde{T}_{n,v}$. For the second see \cite[Proposition 3.9]{SchaeferClassGroups}.
\end{Proof}

Using the Kummer sequence (\ref{Kummer}) one easily determines the size of $\HH^1(k_v,A[n])$. The following lemma is also useful.

\begin{Lemma}
\label{Exponents}
Let $v$ be a nonarchimedean prime, and let $n$ be a positive integer. Then the exponent of $\HH^1(k_v,A[n])$ is equal to the exponent of $A(k_v)/nA(k_v)$.
\end{Lemma}

\begin{Proof}
Let $m$ be the exponent of $A(k_v)/nA(k_v)$ and let $n' = n/m$. By Tate's duality theorem $m$ is also the exponent of $\HH^1(k_v,A)[n]$. Consider the commutative diagram (with exact rows, but non-exact columns)
\[ \xymatrix{ 0 \ar[r]& \frac{A(k_v)}{nA(k_v)} \ar[r]\ar[d]& \HH^1(k_v,A[n])\ar[r]\ar[d]^{m_*} & \HH^1(k_v,A)[n] \ar[r]\ar[d]^{m} & 0 \\
0 \ar[r]& \frac{A(k_v)}{n'A(k_v)} \ar[r]^{\delta_{n'}}\ar[d]^{m}& \HH^1(k_v,A[n'])\ar[r]\ar[d]^{i_*} & \HH^1(k_v,A)[n'] \ar[r]\ar@{^{(}->}[d] & 0\\
0 \ar[r]& \frac{A(k_v)}{nA(k_v)} \ar[r]^{\delta_{n}}& \HH^1(k_v,A[n])\ar[r] & \HH^1(k_v,A)[n] \ar[r]& 0 }\] Here $\delta_n$ and $\delta_{n'}$ denote the connecting homomorphisms in the Kummer sequences for $n$ and $n'$, respectively. The maps $i_*$ and $m_*$ are induced by $A[n'] \stackrel{i}{\hookrightarrow} A[n]$ and $A[n] \stackrel{m}{\to} A[n']$, respectively. The composition of the middle column gives multiplication by $m$ in the group $\HH^1(k_v,A[n])$. Our assumption is that both maps labelled $m$ in the diagram are the zero map. A diagram chase then shows that $i_*\circ m_* = 0$ as well. So $\HH^1(k_v,A[n])$ has exponent dividing $m$.
\end{Proof}

\begin{Proposition}
\label{ImageOb}
Let $v$ be a nonarchimedean prime and $n$ a positive integer. If $A(k_v)/nA(k_v)$ is nontrivial, then $\Delta_{n,v}$ is not the zero map. If the exponent $m$ of $A(k_v)/nA(k_v)$ is odd and the principal polarization is symmetric, then $\Delta_{n,v}\left(\HH^1(k_v,A[n])\right) = \Br(k_v)[m]$.
\end{Proposition}

Applying the Chebotar\"ev density theorem to the extension $k(A[n])/k$ and using lemma \ref{sizetorsors}, the first statement yields the following corollary.

\begin{Corollary}
\label{Deltanot0}
Let $A/k$ be a principally polarized abelian variety over a number field and $n \ge 2$. There are infinitely many primes $v$ of $k$ for which $\Delta_{n,v}$ is not the zero map.
\end{Corollary}

\begin{Proof}
For the first statement recall that the bilinear form associated to $\Delta_{n,v}$ is $T_{n,v}$, which by Tate's theorem is nondegenerate. The assumption implies that $\HH^1(k_v,A[n]) \ne 0$, so the result is clear. For the second, lemma \ref{Exponents} shows that the exponent of $\HH^1(k_v,A[n])$ is equal to that of $A(k_v)/nA(k_v)$. By Zarkhin's result, $\Delta_{n,v}$ is a quadratic form on the finite dimensional $\Z/m\Z$-module $\HH^1(k_v,A[n])$. Since $m$ is odd, the image of $\Delta_{n,v}$ is contained in the $m$-torsion subgroup of $\Br(k_v)\simeq \Q/\Z$ (this is a general property of quadratic maps; see \cite[Lemma 42]{ClarkAV}). We are thus in a position to apply lemma \ref{QuadraticProperties}. We need to show that $\HH^1(k_v,A[n])$ contains an isotropic element of order $m$. The compatibility of $T_{n,v}$ and $\tilde{T}_{n,v}$ shows that the image of any element in $A(k_v)/nA(k_v)$ under the connecting homomorphism is isotropic with respect to $T_{n,v}$. By assumption on the exponent there is such an element of order $m$.
\end{Proof}

The (middle third of) the Poitou-Tate exact sequence for $A[n]$ (see \cite[Section II.6]{GC}) and the fundamental exact sequence for the Brauer group of $k$ give rise to the following commutative diagram with exact rows.
\[ \xymatrix{ 
& \HH^1(k,A[n]) \ar[r]^\res\ar[d]^{\Delta_n} & \prod'\HH^1(k_v,A[n]) \ar[r]^\phi\ar[d]^{\prod \Delta_{n,v}} & \HH^1(k,A[n])^* & \\
0 \ar[r]& \Br(k) \ar[r]^\res& \bigoplus \Br(k_v) \ar[r]^{\sum}& \Q/\Z \ar[r] & 0 }\] The product is the restricted product with respect to unramified subgroups. For all but finitely many primes, the unramified subgroups are exact annihilators with respect to the Tate pairing, so the vertical map in the middle is well-defined.

Clark asked (\cite[Question 45]{ClarkAV}) if $\Delta_n\left(\HH^1(k,A[n])\right) = \Br(k)[n]$ when $A/k$ is a strongly principally polarized abelian variety over a number field. The analysis above shows that this is not the case in general (or even generically; see the remark below). If there is some $\sigma \in G_k$ which acts on $A[n]$ with no nontrivial fixed points, then by Chebotar\"ev's density theorem there are infinitely many primes $v$ for which $A(k_v)/nA(k_v)$, $\HH^1(k_v,A)[n]$ and $\HH^1(k_v,A[n])$ are all trivial (cf. lemma \ref{sizetorsors}). The diagram above shows that every element in the image of $\Delta_n$ must be trivial at these primes. One might still ask the following.

\begin{Question}
\label{ImObQuestion}
Consider an abelian variety $A/k$ over a number field $k$ with a (symmetric) principal polarization. Let $D \in \Br(k)$ be such that the local invariant $\inv_v(D)$ has order dividing the exponent of $A(k_v)/nA(k_v)$ for every prime $v$. Is $D$ in the image of the obstruction map $\Delta_n$?
\end{Question}

It seems unlikely that this is the whole story, since exactness of the Poitou-Tate sequence should place additional restrictions on the image of $\Delta_n$. The results of the next section (e.g. corollary \ref{Density1}) show at least that, given $D$ as in the question, there exist elements in the image of $\Delta_n$ whose local invariants coincide with $D$ on an arbitrarily large finite sets of primes.\\

\begin{Remark} By a result of Serre \cite{Serre} (see also \cite[Corollary 2.1.7]{McQuillen}) there exists an integer $d = d(A/k)$ such that for any $n \ge 1$, all $d$-th powers in $(\Z/n\Z)^\times$ arise as homotheties via the action of $G_k$ on $A[n]$. In particular, for any $n$ sufficiently large, there exists $\sigma \in G_k$ which acts on $A[n]$ with no nontrivial fixed points.
\end{Remark}

We close this section with the following theorem of Lang and Tate \cite[Corollary 1 to Theorem 1, p. 676]{LangTate} which gives a simple criterion for determining splitting fields of torsors over local fields.

\begin{Theorem}[Lang-Tate]
\label{LTtheorem} Let $v$ be a prime of $k$ above the rational prime $q$. Suppose that $A$ has good reduction at $v$, that $n$ is prime to $q$ and that $V \in \HH^1(k_v,A)[n]$ is a torsor under $A/k_v$ of exact period $n$. Then for any finite field extension $K_v/k_v$, $V(K_v) \ne \emptyset$ if and only if $n$ divides the ramification index of $K_v/k_v$.
\end{Theorem}

\section{Weak approximation}
Let $A/k$ be a principally polarized abelian variety over a number field $k$. Under the assumption of full level $n$ structure, one has $A[n] \simeq (\Z/n\Z)^{2g}$ as $G_k$-modules (where $g = \dim A$) and $\mu_n \subset k$. Consequently the cohomology group $\HH^1(k,A[n])$ splits as $2g$ copies of $k^\times/k^{\times n}$. The Grunwald-Wang theorem (or rather a slight variation thereof \cite[Theorems 9.1.11, 9.2.3]{CON}) implies that the restriction map in the following diagram is injective and has dense image in the product of the discrete topologies.
\begin{align}
\label{leveln}
 \xymatrix{ 
\HH^1(k,A[n]) \ar@{=}[r]\ar[d]_{\res}&\left(k^\times/k^{\times n} \right)^{2g} \ar[d]_\res \\
\prod_v\HH^1(k_v,A[n]) \ar@{=}[r]&\prod_v\left(k_v^\times/k_v^{\times n} \right)^{2g} }
\end{align} This, together with the fact that one can interpret the obstruction map $\Delta_n$ in terms of norm residue symbols in this case seems to be responsible for `most' of the period-index discrepancies constructed (see for example, \cite[Proposition 7]{WC1}, \cite{ClSh}, or \cite{ClarkAV}).

Another approach used by Clark (e.g. \cite[Proof of theorem 1]{G1everyindex}) is to make use of the fact \cite[Theorem I.6.26(b)]{ADT} that for any abelian variety $A/k$ over a number field with $\Sha(A/k)[p^\infty]$ finite, there is an exact sequence
\begin{align}
\label{Duality1}
\Sha(A/k)[p^\infty] \to \HH^1(k,A)[p^\infty] \stackrel{\res}{\To} \bigoplus \HH^1(k_v,A)[p^\infty] \to A(k)^{\vee \wedge *}\,,
\end{align} the last term being the Pontryagin dual of the pro-$p$ completion of the Mordell-Weil group of the dual abelian variety. If one is willing to assume that $\Sha(A/k)$ and $A(k)$ are finite, one again gets very satisfactory control over the restriction map. But at present, this limits us to elliptic curves of analytic rank zero over $\Q$. When applied to abelian varieties (i.e. taking $M = A[n]$) the {\em weak weak approximation theorem} below gives a finite level version of (\ref{Duality1}) which requires no assumption on either the Mordell-Weil or Tate-Shafarevich groups.

For the statement we need the following notation. For any $G_k$-module $M$ and any set of primes $S$ of $k$, let $\Sha^1(k,M;S)$ denote the kernel of the restriction map \[\HH^1(k,M) \to \prod_{v\notin S}\HH^1(k_v,M)\,.\] This is the subgroup of cocycle classes that are trivial outside $S$. We use $S^c$ to denote the complement of $S$.  

\begin{Theorem}[weak weak approximation]
\label{weakweakapproximation}
Let $M$ be a finite $G_k$-module and $S$ the finite set of primes where $M^\vee$ is ramified. Then the map 
\[ \res_{S^c}:\Sha^1(k,M;S^c) \to \prod_{v\notin S}\HH^1(k_v,M) \] has dense image in the product of the discrete topologies.
\end{Theorem}

Results of this kind are well known. The statement here is a rather simple generalization of \cite[9.2.2--9.2.3]{CON} or \cite[I.9.8]{ADT}. As we were unable to locate the specific formulation needed for our present application in the literature, we offer the reader \cite[Theorem 1.6]{GW4AV}. We remark that in general it is necessary to exclude the primes in $S$. The (now) standard counterexample led to the correct formulation of the Grunwald-Wang theorem; the map $\HH^1(\Q,\Z/8\Z) \to \HH^1(\Q_2,\Z/8\Z)$ is not surjective. Applying the theorem to $M = A[n]$ and using the criterion of N\'eron-Ogg-Shafarevich one obtains the following corollary (the second result is obtained from the first by use of the Kummer sequence).
\begin{Corollary}
\label{Density1}
Let $A/k$ be an abelian variety over a number field, $n$ an integer and $S$ the finite set of primes containing all primes of bad reduction for $A^\vee$  and all primes dividing $n$. Then the restriction maps
\begin{align*}
\Sha^1(k,A[n];S^c) &\to \prod_{v \notin S}\HH^1(k_v,A[n])\,,\text{ and}\\
\Sha^1(k,A;S^c)[n] &\to \prod_{v \notin S}\HH^1(k_v,A)[n]
\end{align*} have dense images.
\end{Corollary}

\section{Proofs of the main theorems}

It follows from the theorem of Lang and Tate (\ref{LTtheorem}) that any finite collection of torsors of period $n$ that are locally solvable at all `bad primes' can be turned into a finite collection of elements of $\Sha(A/\ell)[n]$ over a sufficiently ramified extension $\ell/k$ of degree $n$. The difficulty is in showing that one can choose the torsors in such a way that they yield distinct, and in particular nontrivial, elements of $\Sha(A/\ell)$. One way of achieving this is to rig the torsors (and their differences) to have index strictly larger than $n$. This is the motivation behind \cite[Theorem 2]{ClSh}, which we generalize to the higher-dimensional situation with the following theorem. This proves theorem \ref{Thm1}.

\begin{Theorem}
\label{PIdifs}
Let $n > 1$ be an integer, $A/k$ a strongly principally polarized abelian variety over a number field $k$ whose N\'eron-Severi group is trivial as a $G_k$-module and $S$ any finite set of primes. There exists a sequence $\{V_i\}_{i=0}^\infty \subset\HH^1(k,A)[n]$ of torsors such that
\begin{enumerate}
\item $V_0(k) \ne \emptyset$.
\item for all $i \ge 0$, $V_i \in \Sha^1(k,A;S^c)[n]$.
\item for all $i \ne j$, $V_i - V_j$ cannot be split by any extension of degree $n$.
\end{enumerate}
If $n$ is a prime power, then the sequence may be chosen so that for all $i \ne j$, the index of $V_i-V_j$ exceeds $n$.
\end{Theorem}

\begin{Proof}
By enlarging $S$ if necessary we may assume $S$ contains all primes of bad reduction for $A$, all primes dividing $n$ and all archimedean primes. Choose representatives $Q_1,\dots,Q_R \in A(k)$ for the finitely many classes in $A(k)/nA(k)$ (Mordell-Weil Theorem), and denote their images in $\HH^1(k,A[n])$ under the connecting homomorphism in the Kummer sequence by $q_1,\dots,q_R$. We will define the sequence inductively. The base of induction is established by setting $V_0 = A$. So assume we have torsors $V_0,\dots,V_{N-1}$ satisfying the three conditions in the theorem, and let $\{ \eta_0,\dots,\eta_{N-1}\} \subset \HH^1(k,A[n])$ denote a set of Kummer lifts of the $V_i$. 

Corollary \ref{Deltanot0} guarantees that there are infinitely many primes $v$ of $k$ for which $\Delta_{n,v}:\HH^1(k_v,A[n]) \to \Br(k_v)$ is not the zero map. We may choose $N\cdot R$ distinct such primes, none of which lie in $S$. Arrange these primes into $N$ sets $\{v_{i,r}\}_{r=1}^R$ of size $R$ (so $0 \le i < N$). For each $v_{i,r}$, let $\xi_{i,r} \in \HH^1(k_{v_{i,r}},A[n])$ be a class such that $\Delta_{n,v_{i,r}}(\xi_{i,r}) \ne 0$. 

By corollary \ref{Density1} there exists a class $\eta_N \in \Sha^1(k,A[n];S^c) \subset \HH^1(k,A[n])$ such that for $0 \le i < N$ and $1 \le r \le R$, we have $\res_{v_{i,r}}(\eta_N) =  \xi_{i,r} + \res_{v_{i,r}}(\eta_{i} - q_r)$. Let $V_N$ be the image of $\eta_N$ in $\HH^1(k,A)[n]$. Then $V_N \in \Sha^1(k,A;S^c)[n]$. Let $i \in \{0,\dots,N-1\}$ and consider $V_N - V_i$. Any Kummer lift of $V_N - V_i$ is of the form $\eta_N - \eta_i + q_r$ for some $r \in \{1,\dots,R\}$. The restriction of this lift at the prime $v_{i,r}$ is equal to $\xi_{i,r}$. Our choice was such that $\Delta_{n,v_{i,r}}(\xi_{i,r}) \ne 0$. Compatibility of the obstruction map with the restriction maps shows that $\Delta_n(\eta_N-\eta_i + q_r) \ne 0$. This is the case for every Kummer lift of $V_N - V_i$, so the result follows from theorem \ref{ObThm}.
\end{Proof}

\begin{Remark}
The reader will note that it is only in the final application of theorem \ref{ObThm} that we make use of the assumptions that the polarization is strong and that the N\'eron-Severi group is trivial as a $G_k$-module.
\end{Remark}

\par\noindent{\sc Proof of Theorems \ref{potentialSha} and \ref{bounds1}:} First note that theorem \ref{bounds1} implies theorem \ref{potentialSha}. We can deduce theorem \ref{bounds1} from theorem \ref{PIdifs} as follows. Let $S$ be the finite set consisting of all primes of bad reduction for $A$, all primes above $p$ and all archimedean primes. Let $T$ be the set of primes outside $S$ where $A(k_v)[p] \ne 0$, and let $\{\ell_j\}$ be a sequence of extensions as in the statement of theorem \ref{bounds1}. Let $N \ge 1$ and let $\{V_i\}_{i=0}^\infty$ be the sequence of torsors given by theorem \ref{PIdifs} (applied with $n = p$). By construction $V_0,\dots,V_N \in \Sha(k,A,S^c)[p]$ and their differences cannot be split by any extension of degree $p$. So for all $j \ge 1$, their images in $\HH^1(\ell_j,A)[p]$ are distinct. It is enough to show that when $j$ is sufficiently large, their images lie in the subgroup $\Sha(A/\ell_j)[p]$.

Let $U$ be the finite set of primes $v$ such that for some $i \le N$, $V_i(k_v) = \emptyset$. If $\ell/k$ is any degree $p$ extension totally ramified at all primes of $U$, then, by the theorem of Lang and Tate (\ref{LTtheorem}), $V_0,\dots,V_N$ represent classes in $\Sha(A/\ell)[p]$. On the other hand, lemma \ref{sizetorsors} shows that $U \subset T$. So by assumption the extension $\ell_j/k$ is totally ramified at all primes of $U$, whenever $j$ is sufficiently large. This completes the proof.
{\hspace*{\fill}\nobreak$\Box$\par\hspace{2mm}}

\par\noindent{\sc Proof of Theorem \ref{bounds2}:}
Let $\ell/k$ be an extension of degree $p$ and let $V/k$ be an element of $\HH^1(k,A)$. Assume $\res_{\ell/k}(V) \in \Sha(A/\ell)[p^\infty]$. In other words, $\res_{\ell/k}(V)$ is an arbitrary element of $\Sha_k(A/\ell)[p^\infty]$. Let $S$ be the set of primes $v$ such that $A$ has bad reduction at $v$ or such that $v$ lies above $p$, and set \begin{align*} T &= \{ v \notin S\,:\, A(k_v)[p] \ne 0 \}\,,\\ U &= S \cup \{ v \in T \,:\, \text{ $v$ is totally ramified in $\ell$ }\}\,. \end{align*}

Since $\res_{\ell/k}(V)$ has $p$-powered period, it also has $p$-powered index. As $\ell/k$ is of degree $p$, we see that $V \in \HH^1(k,A)[p^\infty]$. We claim that the image of $V$ under the map \[\res = \prod_v \res_v :\HH^1(k,A)[p^\infty] \to \bigoplus_v\HH^1(k_v,A)[p^\infty] \] lands in the subgroup \[ G := \prod_{v \in U}\HH^1(k_v,A)[p] \times \prod_{v \notin U}\{ 0 \}\,.\] To prove the claim first note that for all primes $v$, the period and the index of $\res_v(V)$ divide $p$ (since $V$ becomes everywhere locally solvable after an extension of degree $p$). On the other hand, the Lang-Tate result implies that the period and the index of $\res_v(V)$ are equal to $1$ at all nonarchimedean primes $v \notin U$. Now suppose $v$ is archimedean and $w$ is a prime of $\ell$ above $v$. The only nontrivial situation is when $\ell_w \ne k_v$, which can only occur when $v$ is real and ramified, hence in $U$. This establishes the claim.

The bound will be obtained by computing the size of the preimage of $G$ under $\res$. First note that for the archimedean primes in $U$ we have $\#\HH^1(\R,A)[2] = \#\pi_0(A(\R)) \le \#A(\R)[2]$ (see \cite[I.3.7]{ADT}). For nonarchimedean primes the size of $\HH^1(k_v,A)[p]$ is given by lemma \ref{sizetorsors} and theorem \ref{TateLocal}. From this we get
\begin{align*}
\#G &= \prod_{v\in U}\#\HH^1(k_v,A)[p]\\
&=\left(\prod_{v \mid p}\#\HH^1(k_v,A)[p]\right)\cdot\left(\prod_{U \ni v \nmid p}\#\HH^1(k_v,A)[p]\right)\\
&=\left(\prod_{v \mid p}p^{g[k_v:\Q_p]}\#A(k_v)[p]\right)\cdot\left(\prod_{U \ni v \nmid p}\#A(k_v)[p]\right)\\
&= p^{g[k:\Q]}\cdot\prod_{v \in U}\#A(k_v)[p]\\
&\le p^{g([k:\Q] + 2\cdot\#U)} = p^{g[k:\Q] + 2g(\#T + \#S)}\,.
\end{align*} On the other hand, the kernel of $\res$ is the $p$-primary part of $\Sha(A/k)$. This gives the result.
{\hspace*{\fill}\nobreak$\Box$\par\hspace{2mm}}

\section{Examples}
We maintain the notation above. Consider the mod $p$ representation $G_k \to \GL(A[p]) \simeq \GL_{2g}(\F_p)$ associated to $A/k$, where $g=\dim(A)$. If $\sigma_v \in \GL_{2g}(\F_p)$ denotes the image of the Frobenius element at some prime $v \notin S$, then $v \in T$ if and only if $\sigma_v$ has $1$ as an eigenvalue. It is well known (see the remark following question \ref{ImObQuestion}) that if $p$ is sufficiently large, then the image of the mod $p$ representation contains nontrivial scalar matrices. When this is the case, Chebotar\"ev's density theorem gives a positive density set of primes outside of $T$. For elliptic curves without complex multiplication, an earlier result of Serre gives that the mod $p$ representation is surjective for all but finitely many $p$. In this case $T$ has density \[ \frac{ \#\{ g \in \GL_2(\F_p)\,:\, 1 \text{ is an eigenvalue of $g$ } \}}{\#\GL_2(\F_p)} = \frac{(p^2-2)}{(p^2-1)(p-1)}\,,\] which approaches $0$ as $p$ goes to $\infty$. This is of some interest since theorems \ref{bounds1} and \ref{bounds2} are stronger when $T$ is small. On the other extreme, if $A(k)[p] \ne \emptyset$, one has that every prime outside $S$ lies in $T$.

For a concrete example consider the elliptic curve $E/\Q$ of conductor $11$ with minimal model $y^2 + y = x^3 - x^2 - 10x - 20\,.$ For $p = 3$, the bad primes are $S = \{3,11\}$. The mod $3$ representation is surjective so $T$ has density $7/16$. The first few primes in $T$ are $29,$ $53$, $67$, $71$, $79$, $83$, $89,\,\dots\,$. One easily checks that $E(\Q_3)[3] = E(\Q_{11})[3] = 0$, so the group $G$ in the proof of \ref{bounds2} reduces to $\HH^1(\Q_3,E)[3]\times \prod_{q\ne 3}\{0\}$, which has order $3$. A $3$-descent (implemented in MAGMA, see \cite{SchaeferStoll}) shows that $\Sha(E/\Q)[3^\infty] = 0$ (in fact, deeper results prove that $\Sha(E/\Q)= 0$). Thus, there are at most $3$ elements in $\HH^1(\Q,E)$ which restrict into $G$. In fact, there are exactly $3$; the inverse pair of nontrivial trivial torsors are represented by the plane cubic \[ C : -4x^3 + 3x^2z - 6xy^2 - 30xyz - 27xz^2 - 11y^3 - 9y^2z - 12yz^2 + 8z^3 = 0\,,\] which is locally solvable at all primes $q \ne 3$. It follows that if $\ell$ is any cubic number field which is not totally ramified at any primes of $T$, we have \[\#\Sha_\Q(E/\ell)[3^\infty] = \left\{\begin{array}{cc} 3 & \text{ if $C(\ell) = \emptyset$, but $C(\ell_w)\ne\emptyset$ for all $w \mid 3$,}\\ 1 &\text{ otherwise. }\end{array}\right.\] At the same time, theorem \ref{bounds1} gives that $\#\Sha_\Q(E/\ell)[3]$  is unbounded as $\ell$ runs through the fields $\Q(\sqrt[3]{29}),\,$ $\Q(\sqrt[3]{29\cdot 53}),\,\Q(\sqrt[3]{29\cdot 53\cdot 67}),\,\dots\,$.

On the other hand, $E(\Q)$ is generated by the $5$-torsion point $(5,5) \in E(\Q)$. So, for $p = 5$, we have $T = S^c$. For cyclic extensions $\ell/\Q$, Matsuno's result \cite[Corollary 4.5]{Matsuno2} shows that $\dim_{\F_5}\Sha(E/\ell)[5] + \rk (E(\ell))$ grows linearly in the number of ramified primes. Our upper bound for $\dim_{\F_5}\Sha_\Q(E/\ell)[5]$ also grows linearly with the number of totally ramified primes, and one would expect that this is true of $\dim_{\F_5}\Sha_\Q(E/\ell)[5]$ itself. However, theorem \ref{bounds1} only applies to a sequence of extensions which is eventually totally ramified at {\em almost every} prime of $k$. This deficiency ultimately traces back to the use of theorem \ref{weakweakapproximation}, itself a consequence of Poitou-Tate duality which is ineffective as an existence theorem.


\begin{thebibliography}{MM}

\bibitem{Boelling}
  R. B\"olling,  Die Ordnung der Schafarewitsch-Tate Gruppe kann beliebig gro\ss\,werden, Math. Nachr. 67 (1975) 157--179.

\bibitem{Ca1}
 J.W.S. Cassels,  Arithmetic on curves of genus 1, VI. The Tate-\v{S}afarevi\v{c} group can be arbitrarily large, J. Reine Angew. Math. 214--215 (1964) 65--70.

\bibitem{WC1}
  P.L. Clark,  The period-index problems in WC-groups I: elliptic curves, J. Number Theory 114 (2005) 193--208.

\bibitem{G1everyindex}
  P.L.~Clark,  There are genus one curves of every index over every number field, J. Reine Angew. Math. 594 (2006) 201--206.

\bibitem{ClarkAV}
  P.L.~Clark,  Period-index problems in WC-groups II: abelian varieties, (preprint) arXiv:math/0406135v1.
\bibitem{ClSh}
  P.L. Clark and S. Sharif,  Period, index and potential $\Sha$, Algebra and Number Theory 4 (2010) 151--174.

\bibitem{GW4AV}
  B. Creutz,  A Grunwald-Wang type theorem for abelian varieties, (preprint) arXiv:math/1009.3546v2.

\bibitem{Donnelly}
  S. Donnelly,  Elements of given order in Tate-Shafarevich groups of elliptic curves, PhD Thesis (2003) University of Georgia.

\bibitem{Fisher}
  T. Fisher,  Some examples of $5$ and $7$ descent for elliptic curves over $\Q$, J. Eur. Math. Soc. 3 (2001) 169--201.

\bibitem{Kloosterman}
  R. Kloosterman,  The p-part of Tate-Shafarevich groups of elliptic curves can be arbitrarily large, J. Th\'eor. Nomb. Bordeaux 17 (2005) 787--800.

\bibitem{KlSc}
  R. Kloosterman and E.F. Schaefer,  Selmer groups of elliptic curves that can be arbitrarily large, J. Number Theory 99 (2003) 148--163.

\bibitem{Kramer}
  K.~Kramer, A family of semistable elliptic curves with large Tate-Shararevitch groups, Proc. Amer. Math. Soc., 89 (1983) 379--386.

\bibitem{LangTate}
 S. Lang and J. Tate,  Principal homogeneous spaces over abelian varieties, Amer. J. Math. 80 (1958) 659-684.

\bibitem{Matsuno}
  K.~Matsuno,  Construction of elliptic curves with large Iwasawa $\lambda$-invariants and large Tate-Shafarevich groups, Manuscripta Math. 122 (2007) 289--304.

\bibitem{Matsuno2} 
 K.~Matsuno,  Elliptic curves with large Tate-Shafarevich groups over a number field, Math. Res. Lett. 16 (2009) 449--461.

\bibitem{McQuillen}
  M.~McQuillen, Division points on semi-abelian varieties, Invent. Math. 120 (1995) 143-159.

\bibitem{ADT}
   J.S. Milne, Arithmetic duality theorems, Perspectives in Mathematics 1, Academic Press, Boston, 1986.

\bibitem{Mumford}
   D. Mumford,  On the equations defining abelian varieties. I, Invent. Math. 1 (1966) 287--354.

\bibitem{CON}
   J. Neukirch, A. Schmidt and K. Wingberg, Cohomology of number fields (second edition),  Grundlehren Math. Wiss. 323, Springer-Verlag, Berlin, 2008.

\bibitem{ONeil}
 C.H. O'Neil,  The period-index obstruction for elliptic curves, J. Number Theory 95 (2002) 329--339.
 
\bibitem{SchaeferClassGroups}
 E.F. Schaefer,  Class groups and Selmer groups, J. Number Theory 56 (1996) 79--114. 

\bibitem{SchaeferStoll}
  E.F. Schaefer and M. Stoll,  How to do a $p$-descent on an elliptic curve, Trans. Amer. Math. Soc. 356 (2004) 1209-1231.

\bibitem{Serre}
  J-P. Serre,  Quelques propri\'et\'es des groupes alg\'ebriques commutatifs, Appendix in Ast\'erisque 69-70 (1979) 191-202.

\bibitem{GC}
   J-P. Serre, Galois cohomology (second printing), Springer Monogr. Math., Springer-Verlag Berlin Heidelberg New York, 2002.

\bibitem{vanderGeer} 
  G. van der Geer and B. Moonen,  Abelian varieties, available at: \text{http://staff.science.uva.nl/$\sim$bmoonen/boek/BookAV.html}

\bibitem{Zarkhin} 
  Yu.G. Zarkhin,  Noncommutative cohomologies and Mumford groups, Math. Notes 15 (1974) 241--244. Translated from Mat. Zametki 15 (1974) 415--419.

\end{thebibliography}
\end{document}